\input xy \xyoption{all}
\documentclass[12pt,leqno]{amsart}

\usepackage{amscd,amssymb,amsmath,latexsym}

\title{\sc Free and properly discontinuous actions of groups on homotopy $2n$-spheres}

\author{Marek Golasi\'nski, Daciberg Lima Gon\c{c}alves and Rolando Jim\'enez}

\def\R{{\mathbb R}}
\def\S{{\mathbb S}}
\def\Z{{\mathbb Z}}

\newtheorem{thmX}{Theorem}[section]
\newtheorem{proX}[thmX]{Proposition}
\newtheorem{lemX}[thmX]{Lemma}
\newtheorem{CorX}[thmX]{Corollary}
\newtheorem{RemX}[thmX]{Remark}

\newtheorem{questX}[thmX]{Question}

\begin{document}

\vspace*{3cm}

\begin{abstract}
Let $G$ be a group acting freely, properly discontinuously  and cellularly on a
finite dimensional $C$W-complex $\Sigma(2n)$ which has the homotopy type of the $2n$-
sphere $\mathbb{S}^{2n}$. Then, this action induces  an action of the group $G$ on the top cohomology of $\Sigma(2n)$.
For the family  of virtually cyclic groups,  we classify all groups which act on $\Sigma(2n)$, the homotopy type of all possible orbit
spaces and all actions on the top cohomology as well.
\par  Under the hypothesis that  $\mbox{dim}\,\Sigma(2n)\leq 2n+1$, we study the groups
with the virtual cohomological dimension $\mbox{vcd}\,G<\infty$ which act as above on
$\Sigma(2n)$. It turns out that they consist of free groups  and certain semi-direct products $F\rtimes \mathbb{Z}_2$ with $F$ a free group.
For those groups  $G$ and a given  action of $G$ on $\mbox{Aut}(\mathbb{Z})$,  we present an algebraic criterion equivalent to  the realizability
of an action   $G$ on $\Sigma(2n)$ which induces the given action  on its top cohomology.
Then, we obtain a classification of those groups together with actions on the top cohomology of $\Sigma(2n)$.
\end{abstract}

\maketitle

\noindent
\footnote{


The first and second authors gratefully acknowledge a support by Instituto de Matem\'aticas,
UNAM, Oaxaca Branch, where the main part of this paper has been discussed. The authors were
partially supported by CONACyT Grant 98697 and the third author was also supported by  CONACyT Grant 151338.}

\vspace{3mm}

\large

\normalsize

{\large\bf Introduction.} The statement of the spherical space form problem in dimension $n$ is: {\em classify
all manifolds with the $n$-sphere $\mathbb{S}^n$ as the universal cover.} Consequently, manifolds with finite
fundamental groups. The development of that motivates classifications of the possible groups (not necessarily finite) which
act freely, properly discontinuously and cellularly on an $n$-homotopy sphere $\Sigma(n)$ (a finite dimension $CW$-complex with the homotopy type of
the $n$-sphere $\mathbb{S}^n$).  Further, this development began to accelerate with the discovery
by J.\ Milnor \cite{Mi} that some periodic groups could not act freely on any sphere.
Then, R.\ Swan \cite{S0} showed that every periodic finite group acts freely on a finite $CW$-complex homotopic to $\mathbb{S}^{km-1}$
for some $k$, where $m$ is the period of the group. Finite groups with such actions on $\Sigma(n)$ have been fully classified
by Suzuki-Zassenhasus,  see e.g. \cite[Chapter IV, Theorem 6.15]{AM}.
\par A free action of a discrete (finite or infinite) group $G$ on $\Sigma(n)$ induces an action on $H^n(\Sigma(n),\mathbb{Z})$, i.e.,
a homomorphism $G \to \mbox{Aut}(H^{n}(\Sigma(n),\mathbb{Z}))$. Following \cite[Proposition 10.2]{B}, for any action of a
finite group $G$ on $\Sigma(2n+1)$, the induced action $G \to \mbox{Aut}(H^{2n+1}(\Sigma(2n+1),\mathbb{Z}))$ is trivial.
On the other hand, in view of \cite{S1}, the only finite groups acting freely on $\Sigma(2n)$ are, up to isomorphism, trivial or $\mathbb{Z}_2$ and the induced homomorphism $\mathbb{Z}_2\to\mbox{Aut}\,(H^{2n}(\Sigma(2n),\mathbb{Z}))$ is non-trivial.
If the group $G$ is infinite there are more possibilities for the induced action of $G$ on $H^n(\Sigma(n),\mathbb{Z})$ than in the finite case,
and that is a part of the problem to characterize those induced actions.
\par Actions of infinite discrete groups on $\Sigma(n)$ have been also studied, where the induced action on its top cohomology is in general non-trivial and it    is an important part of the data.
We state below some of the relevant results about this study motivated by a problem raised by C.T.C.\ Wall \cite[p.\ 518]{W1}: {\em whether any countable group with periodic Farell cohomology can act freely and properly
on some product $\mathbb{R}^m\times\mathbb{S}^n$?}
\par A breakthrough on Wall's question was made by Connolly and Prassidis (1989). In view of \cite[Corollary 1.4]{C-P}, {\em a discrete group $G$ with the virtual cohomological dimension $\mbox{\em vcd}\,G<\infty$
acts freely and properly on $\mathbb{R}^m\times\mathbb{S}^n$ for some $m,n$ if and only if $G$ is countable and the  Farrell cohomology
$\hat H^\ast(G,-)$ is periodic} (see \cite{C-P} for the definition of the periodicity).
Using arguments presented in \cite{C-P}, it was proved in \cite[Proposition 9.1 and Proposition 9.3]{L} that if a group $G$ which is
not torsion-free with $\mbox{vcd}\,G <\infty$ acts freely and properly discontinuously on $\R^m\times \S^n$ then the period of $\hat H^\ast(G,\Z)$
is two or divides $n + 1$ provided either $n$ is even or odd, respectively.

\par It follows from \cite{AS}  that a discrete group $G$ has {\em periodic cohomology} (after $d$-steps with $d\geq 0)$ if there is
a positive integer $q$ and a class $\alpha\in H^q(G,\mathbb{Z})$
such that the cup product map
$$\alpha\cup-: H^i(G,M)\longrightarrow H^{i+q}(G,M)$$
is an isomorphism for every $G$-module $M$ and $i> d$.
The result \cite[Corollary 1.3]{AS} characterizes groups which act freely and properly discontinuously on $\mathbb{R}^m\times \mathbb{S}^n$,
extends Wall's question above for groups with finite virtual cohomological dimension and states:
{\em A discrete group $G$ acts freely and properly  on $\mathbb{R}^m\times\mathbb{S}^n$ for some $m,n >0$
if and only if $G$ is a countable group with periodic cohomology.}

Further, the result of Johnson \cite[Theorem on p.\ 387]{JJ} states:

\noindent
{\em Let $G$ be a group. Then the following are equivalent:

\mbox{\em (i)} there is a manifold $M$ of type $K(G,1)$;

\mbox{\em (ii)} there is a covering action of $G$ on $\mathbb{R}^m$ for some $m$;

\mbox{\em (iii)} $G$ is countable and has finite cohomological dimension.}

\noindent
Consequently, such a group $G$ acts freely and properly discontinuously on $\mathbb{R}^m\times\mathbb{S}^n$ for any $n >0$.
For more about this  subject, we refer the reader to the papers \cite{AS}, \cite{C-P}, \cite{MT},
\cite{PS} and \cite{W}.

\par By O.\ Talelli \cite{T1}, a group $G$ is said to have {\em periodic cohomology} (after $d$-steps with $d\ge 0$) if there is
a positive integer $q$ such that the functors $H^i(G,-)$ and $H^{i+q}(G,-)$ are naturally equivalent for $i>d$. The class of finite periodic groups
has been extended (see e.g., \cite{MT} and \cite{T2}) to a larger class of discrete groups for which natural equivalences $H^i(G,-)\cong H^{i+q}(G,-)$
of functors for $i>d$ are given by cup product. We point out that it is an open problem
if that periodicity is always induced by cup product with a cohomology class (see e.g., \cite[Remark 2.12]{AS}).
\par It is not clear how to apply most of the results and techniques which appear in \cite{AS} and \cite{C-P} to the cases $n=1$ and $n$ even.
The study of properly discontinuous and cellular actions of discrete groups on a homotopy circle $\Sigma(1)$ was done in \cite{GGJ}
using  different methods than those in the papers mentioned above.

The purpose of this paper is to study free, properly discontinuous and cellular actions of infinite groups $G$ on $\Sigma(2n)$. This study takes also into account the induced actions
of $G$ on $H^{2n}(\Sigma(2n),\Z)$.
\par In virtue of \cite[Proposition 7.1]{L}: {\em the group $G$ is torsion-free or $G\cong G_0\rtimes\mathbb{Z}_2$ with a torsion-free
subgroup $G_0$ provided $G$ acts freely and properly discontinuously on $\mathbb{R}^m\times\mathbb{S}^{2n}$.}

For actions of virtually cyclic groups we show:

{\bf Proposition \ref{Pp}.} {\em Let $G\times \Sigma(2n)\to \Sigma(2n)$ be an action of a non-trivial virtually cyclic
group $G$ on $\Sigma(2n)$ and $\varphi: G \to \mbox{\em Aut}\,(H^{2n}(\Sigma(2n),\mathbb{Z}))$ the induced homomorphism. Then:

\mbox{\em (1)} $G$ is isomorphic to one of the groups: $\mathbb{Z}_2$, $\mathbb{Z}$, $\mathbb{Z}\oplus\mathbb{Z}_2$ or $\mathbb{Z}_2\ast\mathbb{Z}_2\cong\mathbb{Z}\rtimes\mathbb{Z}_2$;

\mbox{\em (2)} any of the groups above admits an action on some $\Sigma(2n)$ and the pair $(G, \varphi)$ is realizable provided:

\hspace{2mm}
\mbox{\em (i)} $G\cong\Z_2$ and $\varphi$ is non-trivial;

\hspace{2mm}
\mbox{\em (ii)}  $G\cong\Z$ and $\varphi$ is any homomorphism;

\hspace{2mm}
\mbox{\em (iii)}  $G\cong\Z\oplus\Z_2$, the restriction $\varphi|_{\Z}$ is trivial
and  $\varphi|_{\Z_2}$ is non-trivial;

\hspace{2mm}
\mbox{\em (iv)} $G\cong\Z\rtimes\Z_2$, the restriction $\varphi|_{\Z}$ is trivial
and  $\varphi|_{\Z_2}$ is non-trivial.
\par Further, the orbit space $\Sigma(2n)/G$ has the homotopy type of one of the manifolds:
$\mathbb{R}P^{2n}$, $\mathbb{S}^{2n}\times\mathbb{S}^1$,
$\mathbb{S}^{2n}\tilde{\times}\mathbb{S}^1$ $($the only non-trivial $\mathbb{S}^{2n}$-bundle over $\mathbb{S}^1$$)$, $\mathbb{R}P^{2n}\times\mathbb{S}^1$
or $\mathbb{R}P^{2n+1}\sharp\mathbb{R}P^{2n+1}$.}

\par Let $F$ be a free group. Given homomorphisms $\theta :\mathbb{Z}_2\to\mbox{Aut}\,(F)$ and
$\varphi : F\rtimes_\theta\mathbb{Z}_2\to\mbox{Aut}\,(\mathbb{Z})\cong\mathbb{Z}_2$ with $\varphi|_{\mathbb{Z}_2}=\mbox{id}_{\mathbb{Z}_2}$,
we say that the pair $(\theta,\varphi)$ is {\em realizable} if there is  an action
$$(F\rtimes_\theta\mathbb{Z}_2)\times\Sigma(2n)\to\Sigma(2n)$$
such that the induced homomorphism $F\rtimes_\theta\mathbb{Z}_2\to\mbox{Aut}\,(H^{2n}(\Sigma(2n),\mathbb{Z}))$
coincides with $\varphi : F\rtimes_\theta\mathbb{Z}_2\to\mbox{Aut}\,(\mathbb{Z})$. The key Lemma \ref{lem} states the necessary and sufficient
conditions for a pair  $(\theta,\varphi)$ to be realizable.
\par For a free group $F_m$ of finite rank $m\ge 1$, we define $m\times m$-matrices $A(k,r,s)$ over
the integers which satisfy $A(k,r,s)^2=I_m$
with $k$ matrices $\begin{pmatrix}0&1\\1&0\end{pmatrix}$, the identity matrix $I_r$ and $-I_s$ on the diagonal for $m=2k+r+s$.
Then, we make use of the the well-known representation $\rho_m : \mbox{Aut}\,(F_m)\to GL_m(\mathbb{Z})$
to prove the following:

{\bf Theorem \ref{p}.} {\em  Let $F_m=\big<x_1,\ldots,x_m\big>$ be  a free group with $m\ge 1$, $\theta
: \mathbb{Z}_2\to \mbox{\em Aut}(F_m)$ and $\varphi : F_m\rtimes_\theta
\mathbb{Z}_2\to\mbox{\em Aut}(\mathbb{Z})$ be homomorphisms  such that $\rho_m(\theta(1_2))=A(k,r,s)$ and
${\varphi}|_{\mathbb{Z}_2}=\mbox{\em id}_{\mathbb{Z}_2}$. Then the pair $(\theta,\varphi)$ is realizable
if and only if $\varphi(x_l,0)=0$ for $l=2k+r+1,\ldots,2k+r+s$.}

\noindent  The result \cite[Theorem 3]{DS} shows that given a free group $F_m$ and $\theta$
it is always possible to find a basis $\{x_1,\ldots,x_m\}$ for $F_m$ which satisfies
the hypothesis of the theorem above.

\par This paper is  organized into four sections, besides this Introduction.

In Section 1 basic facts on actions of groups on $\Sigma(2n)$
are presented.
\par Section 2 aims to determine all non-trivial virtually cyclic groups admitting such actions on $\Sigma(2n)$ and homotopy types
of orbit spaces. The main result is Proposition \ref{Pp}, but Corollary 2.3 classifies orbit spaces of free and properly discontinuous
actions of finite groups on certain spaces having universal covering $\mathbb{R}\times\S^{2n}$.
\par Section 3 analyzes actions $G\times\Sigma(2n)\to\Sigma(2n)$ with the virtual
cohomological dimension $\mbox{vcd}\,G<\infty$ and $\dim\Sigma(2n)\le 2n+1$. We prove
Theorem \ref{p} and make some comments in Remark 3.5 about a version of this theorem for $F$ a free group of infinite rank.
\par Finally, in Section 4 we pose a question  about actions of infinite discrete groups on
$\mathbb{R}^m\times\mathbb{S}^{2n}$.

\par {\bf Acknowledgments.} The authors are indebted to F.X.\ Connolly and S.\ Prasidis for fruitful discussions
on many aspects of this paper, in particular, on the current status of the Wall problem \cite[p.\ 518]{W1} summarized in Section 4.

\renewcommand{\thesection}{{}}
\section{}
\renewcommand{\thesection}{\arabic{section}}

{\large\bf 1.\ Preliminaries.}
A $CW$-complex $\Sigma(n)$ is said to be an $n$-{\em homotopy
sphere}, if $\dim \Sigma(n)<\infty$ and there is a homotopy
equivalence $\Sigma(n)\simeq\mathbb{S}^n$ for the $n$-sphere
$\mathbb{S}^n$ with $n\ge 1$.
\par From now on, we assume that any action $G\times\Sigma(n)\to\Sigma(n)$ of a group $G$ on $\Sigma(2n)$
is free, properly discontinuous and cellular. In the beginning of \cite[Section 1]{GGJ}, we have stated:
\begin{RemX} {\em Notice that $n\le \dim\Sigma(n)$ and for an action
$G\times\Sigma(n)\to\Sigma(n)$ there is a fibration
$$\Sigma(n)\to\Sigma(n)/G\to K(G,1).$$
Consequently, there are isomorphisms
$\pi_k(\Sigma(n))\cong\pi_k(\Sigma(n)/G)$ for $k>1$ and $n\ge 1$,
$\pi_1(\Sigma(n)/G)\cong G$ for $n>1$ and there is an extension
 $$e\to\mathbb{Z}\to\pi_1(\Sigma(1)/G)\to G\to e$$ of groups.\label{R}}
\end{RemX}

Write $\mbox{cd}\,G$ (resp.\  $\mbox{vcd}\,G$) for cohomological (resp.\ virtual
cohomological) dimension of a group $G$ \cite[Chapter VIII]{B}.
\par Given an action $G\times\Sigma(2n)\to\Sigma(2n)$, we consider
the induced homomorphism $$G\to\mbox{Aut}\,(H^{2n}(\Sigma(2n),\mathbb{Z}))\cong\mathbb{Z}_2,$$
which we call from now on the {\em orientation} of the $G$-action.

Then, we make use of \cite[Proposition 7.1]{L} and \cite{S1} to show:
\begin{proX} Let $G\times \Sigma(2n)\to \Sigma(2n)$ be an action of a group $G$ on $\Sigma(2n)$.
Then:\label{pro}

\mbox{\em (1)} $G\cong\mathbb{Z}_2$ or $G=E$ provided $G$ is finite. Further,
$\mathbb{Z}_2\to\mbox{\em Aut}\,(H^{2n}(\Sigma(2n),\mathbb{Z}))$ is

non-trivial;\label{pr}

\mbox{\em (2)} $G$ is torsion-free or $G\cong G_0\rtimes\mathbb{Z}_2$ for some torsion-free subgroup group $G_0$ of $G$.
\end{proX}
{\bf Proof.} (1) If $G$ is finite,  then by \cite{S1},  $G\cong\mathbb{Z}_2$ or $G=E$.
Suppose that $\mathbb{Z}_2\to\mbox{Aut}\,(H^{2n}(\Sigma(2n),\mathbb{Z}))$ is trivial.
Then the Leray-Serre spectral sequence $E_2^{p,q}=H^p(\mathbb{Z}_2,H^q(\Sigma(2n),\mathbb{Z}))$
determined by the fibration $$\Sigma(2n)\to\Sigma(2n)/\mathbb{Z}_2\to K(\mathbb{Z}_2,1)$$
collapses. Hence, the group $H^\ast(\Sigma(2n)/\mathbb{Z}_2,\mathbb{Z})$ does not vanish for infinite many
values of $*$, which
contradicts the fact that
 $\dim \Sigma(2n)/\mathbb{Z}_2<\infty$.

(2) Suppose $G$ is not torsion-free. Then, in view of (1), the induced action
$G\to\mbox{Aut}\,(H^{2n}(\Sigma(2n)),\mathbb{Z})$ is onto and
$G_0=\mbox{ker}\,(G\to\mbox{Aut}\,(H^{2n}(\Sigma(2n),\mathbb{Z}))\cong\mathbb{Z}_2)$ is torsion-free.
Further, the extension $$e\to G_0\to G\to\mathbb{Z}_2\to e$$
splits. Consequently, there is an isomorphism $G\cong G_0\rtimes\mathbb{Z}_2$.

\hfill$\square$

Notice that from Proposition \ref{pro} it follows: if $$\varphi : G\cong G_0\rtimes\mathbb{Z}_2\to\mbox{Aut}\,(H^{2n}(\Sigma(2n),\mathbb{Z}))\cong\mathbb{Z}_2$$
is the induced action then the restriction $\varphi|_{\mathbb{Z}_2}=\mbox{id}_{\mathbb{Z}_2}$.

In \cite[Proposition 1.7]{GGJ}, we have shown:
\begin{proX} If $\mbox{\em vcd}\,G<\infty$ and there is an action $G\times \Sigma(n)\to\Sigma(n)$
then $\mbox{\em vcd}\,G\le \dim\Sigma(n)-n$ for $n\ge 1$. In
particular, $G$ is finite provided $\dim \Sigma(n)=n$.\label{PP}
\end{proX}

Then, we deduce:
\begin{CorX} If $G \times \Sigma(2n) \to \Sigma(2n)$ is an action with $\dim\Sigma(2n) \leq 2n+m$ and $\mbox{\em vcd}\,G<\infty$
then $\mbox{\em cd}\, G\le m$ or $G\cong G_0\rtimes\mathbb{Z}_2$ with $\mbox{\em cd}\,G_0\leq m$. In particular, if $m=1$
then the group $G$ is free or $G\cong F\rtimes\mathbb{Z}_2$ for some free group $F$.\label{c}
\end{CorX}
{\bf Proof.} For an action $G \times \Sigma(2n) \to \Sigma(2n)$ with $\dim\Sigma(2n) \leq 2n+m$ and $\mbox{vcd}\,G<\infty$,
Proposition \ref{PP} yields $\mbox{vcd}\,G\le m$. Then, Proposition \ref{pr} and \cite{Ser1} lead to $\mbox{cd}\, G\le m$ or
$G\cong G_0\rtimes\mathbb{Z}_2$ with $\mbox{cd}\,G_0\leq m$.
\par If $m=1$  then, by means of the above, \cite{St} and \cite{S3}, the group $G$ is free or $G\cong F\rtimes\mathbb{Z}_2$
for some free group $F$.
\hfill$\square$

Now, we show that the family of groups $F\rtimes \mathbb{Z}_2$ for a free group $F$ is closed with respect to free products.
\begin{proX} If $F_i$ are free groups for $i\in I$ then there is an isomorphism
$$\ast_{i\in I}(F_i\rtimes\mathbb{Z}_2)\cong F\rtimes\mathbb{Z}_2$$ for some free group $F$.
\end{proX}
{\bf Proof.} Given $\mathbb{Z}_2=\big<a_i\big>$ for $i\in I$, write $\tilde{F}=\big<x_i|\,i\in I\backslash\{i_0\}\big>$
for the free group generated by the set $\{x_i|\,i\in I\backslash\{i_0\}\}$ for a fixed $i_0\in I$.
Further, consider the homomorphism $\theta: \mathbb{Z}_2=\big<b\big>\to\mbox{Aut}(\tilde{F})$
such that $\theta(b)(x_i)=x_i^{-1}$ for $i\in I\backslash\{i_0\}$. Then, the map
$$\varphi : \tilde{F}\rtimes_\theta\mathbb{Z}_2\longrightarrow \ast_{i\in I}\mathbb{Z}_2$$
given by $\varphi(x_i,0)=a_i\ast a_{i_0}$ for $i\in I\backslash\{i_0\}$, and $\varphi(e,b)=a_{i_0}$
leads to an isomorphism $\tilde{F}\rtimes_\theta\mathbb{Z}_2\stackrel{\cong}{\to}\ast_{i\in I}\mathbb{Z}_2$,
 so the group  $\tilde{F}$ can be regarded as a subgroup of the main group via this  isomorphism.
\par Next, consider the split epimorphism $$p : \ast_{i\in I}(F_i\rtimes\mathbb{Z}_2)\longrightarrow\mathbb{Z}_2,$$
where $p|_{F_{i}\rtimes\mathbb{Z}_2}: F_{i}\rtimes\mathbb{Z}_2\to \mathbb{Z}_2$ is the projection map
for all $i\in I$.

\par Notice that $\tilde{F}\ast(\ast_{i\in I}(F_i\times\{0\}))\subseteq \mbox{Ker}\,p$ and
$\tilde{F}\ast(\ast_{i\in I}(F_i\times\{0\}))$ is a normal subgroup of $(\ast_{i\in I}(F_i\rtimes\mathbb{Z}_2)$.
Further, $(g_{i_1},1)\ast\cdots\ast(g_{i_n},1)\in \mbox{Ker}\,p$ yields that $n$ is even.
This shows that $\mbox{Ker}\,p=\tilde{F}\ast(\ast_{i\in I}(F_i\times\{0\})$ and the proof
is complete.

\hfill$\square$

\renewcommand{\thesection}{{}}
\section{}
\renewcommand{\thesection}{\arabic{section}}

{\large\bf 2.\ Virtually cyclic groups acting on $\Sigma(2n)$.}
Recall that a {\em virtually cyclic} group is a group that has a cyclic subgroup of
finite index. The following criterion is mainly due to P.\ Scott and C.T.C.\  Wall \cite{SW}:
 \begin{thmX} Let $G$ be a finitely generated group.
Then, the following are equivalent:

\mbox{\em (1)} $G$ is a group with two ends;

\mbox{\em (2)} $G$ has an infinite cyclic group of finite index;

\mbox{\em (3)} $G$ has a finite normal subgroup $F\trianglelefteq G$
with the quotient $G/F\cong\mathbb{Z}$ or

$\mathbb{Z}_2\star\mathbb{Z}_2\cong D_\infty$, the infinite dihedral group.

\noindent
Equivalently, $G$ is of the form:

\mbox{\em (1)} a semi-direct product $F\rtimes\mathbb{Z}$ with $F$ finite

or\label{T}

\mbox{\em (2)} $G_1\star_F G_2$ with $F$ finite, where $[G_i:F]=2$ for $i=1,2$.
\end{thmX}

\smallskip

\par Given an action $G\times\Sigma(n)\to \Sigma(n)$, write $\alpha_{\Sigma(n)}$
for the first Postnikov invariant \cite{R} of the orbit space $\Sigma(n)/G$.
In the sequel we need:
\begin{lemX} Let a discrete group $G$ act on $\Sigma_1(n)$
and $\Sigma_2(n)$ with $\dim \Sigma_1(n)/G\le n+1$ for $n\ge 2$, and $\dim \Sigma_2(n)/G$ arbitrary.
\par The orbit spaces $\Sigma_1(n)/G$ and $\Sigma_2(n)/G$ have the same homotopy
type if and only if there is an automorphism $\varphi\in\mbox{\em Aut}(G)$
with $\varphi^\ast(\alpha_{\Sigma_2(n)})=\alpha_{\Sigma_1(n)}$.  \label{l}
\end{lemX}
{\bf Proof.} If the orbit spaces $\Sigma_1(n)/G$ and $\Sigma_2(n)/G$ have the same homotopy
type, then certainly there is $\varphi\in\mbox{Aut}(G)$
with $\varphi^\ast(\alpha_{\Sigma_2(n)})=\alpha_{\Sigma_1(n)}$.
\par Now, suppose that there is $\varphi\in\mbox{Aut}(G)$
with $\varphi^\ast(\alpha_{\Sigma_2(n)})=\alpha_{\Sigma_1(n)}$.
Then we derive a map of two stage Postnikov towers of $\Sigma_1(n)/G$ and $\Sigma_2(n)/G$:
$$\xymatrix{
\Sigma_1(n)/G\ar[dr]\ar[r]&X_2\ar[d]\ar[r]^{f_2}&Y_2\ar[d]&\Sigma_2(n)/G\ar[dl]\ar[l]\\
&X_1=K(G,1)\ar[r]^{\bar\varphi}& Y_1=K(G,1),&}$$

\noindent where $\bar \varphi$ is the induced map by $\varphi\in\mbox{Aut}(G)$ at the Eilenberg-MacLane space $K(G,1)$.

Because $\dim \Sigma_1(n)\le n+1$,  obstruction theory leads to a map $f : \Sigma_1(n)/G\to \Sigma_2(n)/G$
with $\pi_n(f)=\pi_n(f_2) :\pi_n(\Sigma_1(n)/G)\to\pi_n(\Sigma_2(n)/G)$ being an isomorphism.
Then, for the lifting $\tilde f$
$$\xymatrix{
\Sigma_1(n)\ar[dd]\ar[rr]^{\tilde{f}}&&\Sigma_2(n)\ar[dd]\\
\\
\Sigma_1(n)/G\ar[rr]^f&& \Sigma_2(n)/G}$$
of $f : \Sigma_1(n)/G\to \Sigma_2(n)/G$,  we deduce that $\pi_n(\tilde{f}):\pi_n(\Sigma_1(n))\to\pi_n(\Sigma_2(n))$ is
an isomorphism. So we have an isomorphism in homology in all dimensions since   the spaces have the homotopy type of a sphere.
Therefore  $\pi_k(\tilde f)$
is an isomorphism for all $k$  and consequently the same for  $\pi_k( f)$. Consequently,
$f : \Sigma_1(n)/G\to \Sigma_2(n)/G$ is a homotopy equivalence,  which completes the proof.

\hfill$\square$

\par Let $G$ be a group and $\varphi : G\to\mbox{Aut}\,(\mathbb{Z})$ a homomorphism. We say that the pair $(G,\varphi)$ is {\em realizable} if there is  an action $G\times\Sigma(2n)\to\Sigma(2n)$
that the induced homomorphism $G\to\mbox{Aut}\,(H^{2n}(\Sigma(2n),\mathbb{Z}))$ coincides with $\varphi : G\to\mbox{Aut}\,(\mathbb{Z})$.
\par Then, we are in position to show:
\begin{proX} Let $G\times \Sigma(2n)\to \Sigma(2n)$ be an action of a non-trivial virtually cyclic
group $G$ on $\Sigma(2n)$ and $\varphi: G \to \mbox{\em Aut}\,(H^{2n}(\Sigma(2n),\mathbb{Z}))$ the induced homomorphism. Then:\label{Pp}

\mbox{\em (1)} $G$ is isomorphic to one of the groups: $\mathbb{Z}_2$, $\mathbb{Z}$, $\mathbb{Z}\oplus\mathbb{Z}_2$ or $\mathbb{Z}_2\ast\mathbb{Z}_2\cong\mathbb{Z}\rtimes\mathbb{Z}_2$;

\mbox{\em (2)} any of the groups above admits an action on some $\Sigma(2n)$ and the pair $(G, \varphi)$ is realizable provided:

\hspace{2mm}
\mbox{\em (i)} $G\cong\Z_2$ and $\varphi$ is non-trivial;

\hspace{2mm}
\mbox{\em (ii)}  $G\cong\Z$ and $\varphi$ is any homomorphism;

\hspace{2mm}
\mbox{\em (iii)}  $G\cong\Z\oplus\Z_2$, the restriction $\varphi|_{\Z}$ is trivial
and  $\varphi|_{\Z_2}$ is non-trivial;

\hspace{2mm}
\mbox{\em (iv)} $G\cong\Z\rtimes\Z_2$, the restriction $\varphi|_{\Z}$ is trivial
and  $\varphi|_{\Z_2}$ is non-trivial.
\par Further, the orbit space $\Sigma(2n)/G$ has the homotopy type of one of the manifolds:
$\mathbb{R}P^{2n}$, $\mathbb{S}^1\times\mathbb{S}^{2n}$,
$\mathbb{S}^1\tilde{\times}\mathbb{S}^{2n}$ $($the only non-trivial $\mathbb{S}^{2n}$-bundle over $\mathbb{S}^1$
being the mapping torus of the antipodal map $\alpha : \S^{2n}\to\S^{2n}$$)$, $\mathbb{S}^1\times\mathbb{R}P^{2n}$
or $\mathbb{R}P^{2n+1}\sharp\mathbb{R}P^{2n+1}$.
\end{proX}
{\bf Proof.} (1): Follows immediately from Proposition \ref{pro} and Theorem \ref{T}.
\par (2):  Writing $\tilde{\mathbb{Z}}$ for the $G$-module structure on $H^{2n}(\Sigma(2n),\mathbb{Z})\cong\mathbb{Z}$
for $G$ being one of the groups from (1), we make use of Lemma \ref{l}.

\par (i):  $G\cong\mathbb{Z}_2$. Then, certainly there is the standard action $\Z_2\times\S^{2n}\to\S^{2n}$ for  any $n\ge 1$ and
by e.g., \cite[Lemma 2.5]{BK}, it holds that  $\Sigma(2n)/\Z_2\simeq\mathbb{R}P^{2n}$,  for any action $\Z_2\times\Sigma(2n)\to\Sigma(2n)$.

\par (ii): $G\cong\mathbb{Z}$. Because $H^{2n+1}(\mathbb{Z},\tilde{\mathbb{Z}})=H^{2n+1}(\mathbb{Z},\mathbb{Z})=0$,
there is at most one homotopy type of $\Sigma(2n)/\Z$ for the non-trivial and trivial $G$-actions on $H^{2n}(\Sigma(2n),\mathbb{Z})$.
Any of them may be realized. Namely, consider the $\mathbb{Z}$-actions:
$$\circ,\bar{\circ} :\mathbb{Z}\times (\mathbb{R}\times\mathbb{S}^{2n})\to \mathbb{R}\times\mathbb{S}^{2n}$$
given by $n\circ(t,x)=(t+n,x)$ and $n\bar\circ(t,x)=(t+n,-x)$, respectively for $n\in\mathbb{Z}$ and $(t,x)\in\mathbb{R}\times\mathbb{S}^{2n}$.
The corresponding orbit spaces are homotopic to
$\mathbb{S}^1\times\mathbb{S}^{2n}$ or $\mathbb{S}^1\tilde{\times}\mathbb{S}^{2n}$, respectively.

\par (iii):  $G\cong\mathbb{Z}\oplus\mathbb{Z}_2$. Then, by Proposition \ref{pro} and its proof,
$G$ acts non-trivially on $H^{2n}(\Sigma(2n),\mathbb{Z})\cong\mathbb{Z}$. Hence, we have an epimorphism $\varphi : \mathbb{Z}\oplus\mathbb{Z}_2\to\mathbb{Z}_2$  such that $\varphi(0,1_2)=1_2$ and  $\varphi(1,0)=0$ or $\varphi(1,0)=1_2$. Further, observe that there is an isomorphism $\mathbb{Z}\oplus\mathbb{Z}_2\stackrel{\cong}{\to}\mathbb{Z}\oplus\mathbb{Z}_2$ given by: $(1,0)\mapsto(1,1_2)$ and $(0,1_2)\mapsto(0,1_2)$.
\par Analysing the Lyndon-Hochschild-Serre spectral sequence corresponding to the extension $$e\to\mathbb{Z}\to\mathbb{Z}\oplus\mathbb{Z}_2\to\mathbb{Z}_2\to e,$$
we deduce that $H^0(\mathbb{Z}\oplus\mathbb{Z}_2,\tilde{\mathbb{Z}})=0$ and $H^k(\mathbb{Z}\oplus\mathbb{Z}_2,\tilde{\mathbb{Z}})\cong\mathbb{Z}_2$ for $k>0$.
In particular, $H^{2n+1}(\mathbb{Z}\oplus\mathbb{Z}_2,\tilde{\mathbb{Z}})\cong\mathbb{Z}_2$
and there are two possible values for the first Postnikov invariant
$K(\mathbb{Z}\oplus\mathbb{Z}_2,1)\to\widehat K(\mathbb{Z},2n+1)$ of the orbit space $\Sigma(2n)/\Z\oplus\Z_2$, where
$\widehat K(\mathbb{Z},2n+1)=\widetilde{K(\mathbb{Z}\oplus\mathbb{Z}_2,1)}\times_{\mathbb{Z}\oplus\mathbb{Z}_2}K(\mathbb{Z},2n+1)$ is
the twisted Eilenberg-MacLane space.
But, the Leray-Serre spectral sequence corresponding to the fibration $K(\mathbb{Z},2n+1)\to X_2\to K(\mathbb{Z}\oplus\mathbb{Z}_2,1)$ shows that this invariant cannot be trivial.
Consequently, there is only one homotopy type of the quotient space realized by the action:
$$\circ :(\mathbb{Z}\oplus\mathbb{Z}_2)\times (\mathbb{R}\times\mathbb{S}^{2n})\to \mathbb{R}\times\mathbb{S}^{2n}$$
given by $(1,0)\circ(t,x)=(t+1,x)$ and $(0,1_2)\circ(t,x)=(t,-x)$ for $(t,x)\in\mathbb{R}\times\mathbb{S}^{2n}$
with the corresponding quotient space homotopic to $\mathbb{S}^1\times\mathbb{R}P^{2n}$.
\par (iv):  $G\cong\mathbb{Z}\rtimes\mathbb{Z}_2$. Then, again by \cite[Proposition 7.1]{L} and its proof,
$G$ acts non-trivially on $H^{2n}(\Sigma(2n),\mathbb{Z})\cong\mathbb{Z}$. Hence, we obtain an
extension $$0\to\mathbb{Z}\to\mathbb{Z}\rtimes\mathbb{Z}_2\to\mathbb{Z}_2\to 0.$$ Then, the corresponding
Lyndon-Hochschild-Serre spectral sequence yields  $H^{2n+1}(\mathbb{Z}\rtimes\mathbb{Z}_2,\tilde{\mathbb{Z}})\cong\mathbb{Z}_2$.
The methods parallel to those in (iii) show that the first Postnikov invariant of the orbit space $\Sigma(2n)/\Z\rtimes\Z_2$
cannot be trivial. Finally, there is also only one homotopy type of the quotient space realized by the action:
$$\circ :(\mathbb{Z}\rtimes\mathbb{Z}_2)\times(\mathbb{R}\times\mathbb{S}^{2n})\to \mathbb{R}\times\mathbb{S}^{2n}$$
given by $(1,0)\circ(t,x)=(t+1,x)$ and $(0,1_2)\circ(t,x)=(-t,-x)$ for $(t,x)\in \mathbb{R}\times\mathbb{S}^{2n}$
with the corresponding orbit space homotopic to $\mathbb{R}P^{2n+1}\sharp\mathbb{R}P^{2n+1}$.
\par The last statement follows from the proof of (2).

\hfill$\square$

In view of \cite[Corollary 2]{T}, the classification of all free and properly discontinuous actions by a finite group on $\S^1\times\S^2$ follows from the observation
that there exist only four compact $3$-manifolds which have $\R\times\S^2$ as a universal covering space.
\par Now, we deduce below that any manifold with the universal covering space $\R\times\mathbb{S}^{2n}$
has the homotopy type one of the following manifolds: $\mathbb{S}^1\times\mathbb{S}^{2n}$, $\mathbb{S}^1\tilde{\times}\mathbb{S}^{2n}$,
$\mathbb{S}^1\times\mathbb{R}P^{2n}$ or $\mathbb{R}P^{2n+1}\sharp\mathbb{R}P^{2n+1}$.
\begin{CorX} Suppose that a finite non-trivial group $G$ acts freely on one of the manifolds:
$\mathbb{S}^1\times\mathbb{S}^{2n}$, $\mathbb{S}^1\tilde{\times}\mathbb{S}^{2n}$, $\mathbb{S}^1\times\mathbb{R}P^{2n}$ or $\mathbb{R}P^{2n+1}\sharp\mathbb{R}P^{2n+1}$. Let $M$ be the orbit
space of $G$.

\noindent
\mbox{\em (1)} If $G$ acts on $\mathbb{S}^1\times\mathbb{S}^{2n}$ then:

\mbox{\em (i)} $G\cong\mathbb{Z}_2$ and $M\simeq\mathbb{S}^1 \times\mathbb{S}^{2n}$, $M\simeq\mathbb{S}^1\tilde{\times}\mathbb{S}^{2n}$,
$M\simeq\mathbb{S}^1\times\mathbb{R}P^{2n}$ or

$M\simeq\mathbb{R}P^{2n+1}\sharp\mathbb{R}P^{2n+1}$;

\mbox{\em (ii)} $G\cong\mathbb{Z}_2\oplus\mathbb{Z}_2$ and $M\simeq\mathbb{S}^1\times\mathbb{R}P^{2n}$ or $M\simeq\mathbb{R}P^{2n+1}\sharp\mathbb{R}P^{2n+1}$;

\mbox{\em (iii)} $G\cong\mathbb{Z}_{2k+1}$ for some $k\ge 1$ and $M\simeq\mathbb{S}^1\times\mathbb{S}^{2n}$;

\mbox{\em (iv)} $G\cong\mathbb{Z}_{2k}$ for some $k>1$  and $M\simeq\mathbb{S}^1\times\mathbb{S}^{2n}$,
$M\simeq\mathbb{S}^1\tilde\times\mathbb{S}^{2n}$ or $M\simeq\mathbb{S}^1\times\mathbb{R}P^{2n}$;

\mbox{\em (v)} $G\cong\mathbb{Z}_2\oplus\mathbb{Z}_{2k}$ for some $k\ge 1$ and $M\simeq\mathbb{S}^1\times\mathbb{R}P^{2n}$;

\mbox{\em (vi)} $G\cong\mathbb{Z}_k\rtimes\mathbb{Z}_2=D_k$ for some $k>2$, the dihedral group of order $2k$ and

$M\simeq\mathbb{R}P^{2n+1}\sharp\mathbb{R}P^{2n+1}$;

\noindent
\mbox{\em (2)} If $G$ acts on $\mathbb{S}^1\tilde{\times}\mathbb{S}^{2n}$ then:

\mbox{\em (i)} $G\cong\mathbb{Z}_{2k+1}$ for some $k\ge 1$ and $M\simeq\mathbb{S}^1\tilde\times\mathbb{S}^{2n}$;

\mbox{\em (ii)} $G\cong\mathbb{Z}_{2k}$ for some $k\ge 1$ and $M\simeq\mathbb{S}^1\times\mathbb{R}P^{2n}$.

\noindent
\mbox{\em (3)} If $G$ acts on $\mathbb{S}^1\times\mathbb{R}P^{2n}$ then
 $G\cong\mathbb{Z}_k$ for some $k\ge 2$ and $M\simeq\mathbb{R}P^{2n}\times\mathbb{S}^1$.

\noindent
\mbox{\em (4)} If $G$ acts on $\mathbb{R}P^{2n+1}\sharp\mathbb{R}P^{2n+1}$ then
$G\cong\mathbb{Z}_2$ and $M\simeq\mathbb{R}P^{2n+1}\sharp\mathbb{R}P^{2n+1}$.

\noindent
Further, in all four cases above, the groups described act on the corresponding manifold.
\end{CorX}
{\bf Proof.} We remark that $\mathbb{R}\times\mathbb{S}^{2n}$ is the universal covering space of
the manifolds listed above and make use of Proposition \ref{Pp}.
\par (1): If $G$ acts on $\mathbb{S}^1\times\mathbb{S}^{2n}$ then the quotient map
$\mathbb{S}^1\times\mathbb{S}^{2n}\to\mathbb{S}^1\times\mathbb{S}^{2n}/G$
is covering and $\mathbb{S}^1\times\mathbb{S}^{2n}/G$ is homotopic to one of those manifolds.
Then, we get an extension of groups $$e\to\mathbb{Z}\to\pi\to G\to e,$$
where $\pi=\pi_1(\mathbb{S}^1\times/G\mathbb{S}^{2n})$.
Certainly, $\mathbb{S}^1\times\mathbb{S}^{2n}/G$ cannot be homeomorphic to $\mathbb{R}P^n$
because the group $G$ is finite.
\par If $\mathbb{S}^1\times\mathbb{S}^{2n}/G\cong\mathbb{S}^1\times\mathbb{S}^{2n}$
or $\mathbb{S}^1\times\mathbb{S}^{2n}\/G\cong\mathbb{S}^1\tilde\times\mathbb{S}^{2n}$
then $G\cong\mathbb{Z}_k$ for some $k\ge 2$.
\par If $\mathbb{S}^1\times\mathbb{S}^{2n}/G\cong\mathbb{S}^1\times\mathbb{R}P^n$
then $G\cong\mathbb{Z}_k$ (for $k$ odd) or $G\cong\mathbb{Z}_2\oplus\mathbb{Z}_k$ (for $k$ even).
\par Because the manifold $\mathbb{S}^1\times\mathbb{S}^{2n}$ is oriented and $(1_2,2k)$ is in the kernel of the epimorphism
$\mathbb{Z}_2\times\mathbb{Z}\to\mathbb{Z}_{4k}\to e$, the group $\mathbb{Z}_{4k}$ cannot act freely on $\mathbb{S}^1\times\mathbb{S}^{2n}$
to obtain  $\mathbb{S}^1\times\mathbb{R}P^{2n}$.
\par If $M\simeq\mathbb{R}P^{2n+1}\sharp\mathbb{R}P^{2n+1}$ then $G\cong\mathbb{Z}_k\rtimes\mathbb{Z}_2$,
the dihedral group of order $2k$ for some $k>2$.

(2): Because $\mathbb{S}^1\tilde\times\mathbb{S}^{2n}$ is non-oriented, $\mathbb{S}^1\tilde\times\mathbb{S}^{2n}/G
\cong\mathbb{S}^1\tilde\times\mathbb{S}^{2n}$
or $\mathbb{S}^1\tilde\times\mathbb{S}^{2n}/G\cong\mathbb{S}^1\times\mathbb{R}P^{2n}$.
Hence, $G\cong\mathbb{Z}_k$ for some $k>1$ or $G\cong\mathbb{Z}_2\oplus\mathbb{Z}_k$
for $k$ odd. But $\mathbb{S}^1\tilde\times\mathbb{S}^{2n}$ is non-oriented, so for the epimorphism
$\mathbb{Z}_2\oplus\mathbb{Z}\to\mathbb{Z}_2\oplus\mathbb{Z}_k$ some element $(1_2,l)$ with $l\not=0$
must be in its kernel. Hence, $G\cong\mathbb{Z}_k$ for some $k>1$, only.

(3): Because $\mathbb{S}^1\times\mathbb{R}P^{2n}$ is non-oriented, $\mathbb{S}^1\times\mathbb{R}P^{2n}/G\cong\mathbb{S}^1\tilde\times\mathbb{S}^{2n}$
or $\mathbb{S}^1\times\mathbb{R}P^{2n}/G\cong\mathbb{S}^1\times\mathbb{R}P^{2n}$.
Then, we obtain that $G\cong\mathbb{Z}_k$ for some $k\ge 2$.

(4): Because $\pi_1(\mathbb{R}P^{2n+1}\sharp\mathbb{R}P^{2n+1})\cong\mathbb{Z}\rtimes\mathbb{Z}_2$, we have
only to analyze the extension $$e\to\mathbb{Z}\rtimes\mathbb{Z}_2\to\mathbb{Z}\rtimes\mathbb{Z}_2\to G\to e.$$
 Because $\mathbb{Z}\rtimes\mathbb{Z}_2$ must be sent to its normal subgroup by the monomorphism $e\to\mathbb{Z}\rtimes\mathbb{Z}_2\to\mathbb{Z}\rtimes\mathbb{Z}_2$, we deduce that $G\cong\mathbb{Z}_2$.

\hfill$\square$

Let $n\ge 2$ and let $\tau$ be a free involution on $\S^1\times \S^n$. Then, in view of \cite[Theorem 2.1]{BK} the quotient $\S^1\times \S^n/\tau$
belongs to one of the four homotopy types: $\S^1\times\S^n$, $\S^1\tilde{\times}\S^n$, $\S^1\times \R P^n$ and $\R P^{n+1}\sharp\R P^{n+1}$
realized by the standard involutions.
\par Now, we are in position to conclude the following generalization of the above, provided $n$ is even:
\begin{CorX}
Let $n\ge 1$ and $\tau$ be a free involution on one the four manifolds:
$\mathbb{S}^1\times\mathbb{S}^{2n}$, $\mathbb{S}^1\tilde{\times}\mathbb{S}^{2n}$, $\mathbb{S}^1\times\mathbb{R}P^{2n}$ or $\mathbb{R}P^{2n+1}\sharp\mathbb{R}P^{2n+1}$.
Then, the corresponding orbit space also belongs  to one of their homotopy types.
\end{CorX}

\renewcommand{\thesection}{{}}
\section{}
\renewcommand{\thesection}{\arabic{section}}

{\large\bf 3.\ Other groups acting on $\Sigma(2n)$.}
Here, we analyze actions $G\times\Sigma(2n)\to\Sigma(2n)$ with  $\mbox{vcd}\,G<\infty$ and
$\dim\Sigma(2n)\le 2n+1$. By  Corollary \ref{c},  the group $G$ is free or
$G\cong F\rtimes\mathbb{Z}_2$ for some free group $F$ with an arbitrary rank.
\par Let $F$ be a free group and $\mbox{Aut}\,(F)$ its automorphism
group. For a homomorphism $\theta :\mathbb{Z}_2\to\mbox{Aut}\,(F)$ we consider the semidirect product
$G\cong F\rtimes_{\theta}\mathbb{Z}_2$ which is completely determined by $\theta$.
Given also $\varphi : F\rtimes_\theta\mathbb{Z}_2\to\mbox{Aut}\,(\mathbb{Z})\cong\mathbb{Z}_2$ with
$\varphi|_{\mathbb{Z}_2}=\mbox{id}_{\mathbb{Z}_2}$, we say that
the pair $(\theta,\varphi)$ is {\em realizable} if $( F\rtimes_\theta\mathbb{Z}_2, \varphi)$ is  realizable
(see Section 1).
\par Notice that any free group $F$ acts on the homotopy $2n$-sphere $(\widetilde{\bigvee_{i\in I}\S^1})\times\S^{2n}$ for any $n\ge 1$,
where $\widetilde{\bigvee_{i\in I}\S^1}$ is the universal covering of the wedge $\bigvee_{i\in I}\S^1$ provided $F=\big<x_i;\ i\in I\big>$.
Consequently, for the trivial homomorphism $\theta_0 :\mathbb{Z}_2\to\mbox{Aut}\,(F)$,
any pair $(\theta_0,\varphi)$ is realizable by the action
$$\circ : (F\rtimes_\theta\mathbb{Z}_2)\times((\widetilde{\bigvee_{i\in I}\S^1})\times\mathbb{S}^{2n})
\to(\widetilde{\bigvee_{i\in I}\S^1})\times\mathbb{S}^{2n}$$
given by: $(g,0)\circ(t,s)=(gt,\mbox{sgn}(g)s,)$  and $(g,1_2)\circ(t,s)=(gt,-\mbox{sgn}(g)s)$ for $g\in F$ and
$(t,s)\in(\widetilde{\bigvee_{i\in I}\S^1})\times\mathbb{S}^{2n}$, where $\mbox{sgn} : F\to\Z_2=\{\pm 1\}$ is the homomorphism determined by
the restriction of $\varphi : F\rtimes_\theta\mathbb{Z}_2\to\mbox{Aut}\,(\mathbb{Z})\cong\mathbb{Z}_2$ to the group $F$.
\par Writing $\mathbb{Z}_2=\big<1_2\big>$, we show a general fact:
\begin{lemX}\mbox{\em (Fundamental Lemma)} The pair $(\theta,\varphi)$ is realizable if and only if it
does not exist  $g\in F$ such that \label{lem}
$$
\left\{\begin{array}{l}
\theta(1_2)(g)=g^{-1},\\
\varphi(g,0)=1_2.
\end{array}
\right.
$$
\end{lemX}
{\bf Proof.}  Let $\theta :\mathbb{Z}_2\to\mbox{Aut}\,(F)$. By the $1$-dimensional analog of the Nielsen realization problem
\cite[Theorems 2.1 and 4.1]{Cul}, the automorphism $\theta(1_2)\in\mbox{Aut}(F)$ can be realized by a homeomorphism
$h : \Gamma\to\Gamma$ of a graph $\Gamma$ with the fundamental group $\pi_1(\Gamma)\cong F$, such that $h$ has a fixed point
and $h^2=\mbox{id}_\Gamma$. Writing $\tilde{\Gamma}$ for the universal covering of $\Gamma$, we get the induced
homeomorphism $\tilde{h} : \tilde{\Gamma}\to\tilde{\Gamma}$ with $\tilde{h}^2=\mbox{id}_{\tilde{\Gamma}}$.
Then, we are in a position to consider a map
$$\circ : (F\rtimes_\theta\mathbb{Z}_2)\times(\tilde{\Gamma}\times\mathbb{S}^{2n})
\to\tilde{\Gamma}\times\mathbb{S}^{2n}$$
given by: $(g,0)\circ(t,s)=((\theta(1_2)g)t,\mbox{sgn}(g)s)$ and $(g,1_2)\circ(t,s)=((\theta(1_2)g)(\tilde{h}(t)),-\mbox{sgn}(g)s)$
for $g\in F$ and $(t,s)\in\tilde{\Gamma}\times\mathbb{S}^{2n}$.
Now, we prove that the map defined above is an action of the group $F\rtimes_\theta\mathbb{Z}_2$.
So for any two elements $w_1,w_2\in F\rtimes_\theta\mathbb{Z}_2$ and
$(t,s)\in \tilde{\Gamma}\times\mathbb{S}^{2n}$ we must show that $w_2(w_1(t,s))=(w_2 w_1)(t,s)$.
\noindent
Notice that $\mbox{sgn}(g\theta(1_2)g')=\mbox{sgn}(gg')$ and
$(g,1_2)(s,t)=((g,0)(e,1_2))(t,s)=(g,0)(\tilde{h}(t),-s)=((\theta(1_2)g)\tilde{h}(t),-\mbox{sgn}(g)s)$.
The case where $w_i=(g_i,\bar 0)$ for $i=1,2$ is easy and we leave to the reader. For the
remaining  cases we have:

\noindent
(i) $((g,0)(g',1_2))(t,s)=(gg',1_2)(t,s)=((\theta(1_2)(gg'))\tilde{h}(t),-\mbox{sgn}(gg')s)$
and

\noindent
$(g,0)((g',1_2))(t,s))=(g,0)=((g,0)(g',1_2))(t,s)=(g,0)((\theta(1_2)g')\tilde{h}(t),-\mbox{sgn}(g')s)=
((\theta(1_2)g)(\theta(1_2)g')\tilde{h}(t),-\mbox{sgn}(g)\mbox{sgn}(g')s)=(g,0)(g',1_2))(t,s)$;

\vspace{1mm}

\noindent
(ii) $((g,1_2)(g',0))(t,s)=(g\theta(1_2)g',1_2)(t,s)=(\theta(1_2)(g\theta(1_2)g')\tilde{h}(t),-\mbox{sgn}(g\theta(1_2)g')s)=
((\theta(1_2)g)g'\tilde{h}(t),\mbox{sgn}(gg')s)$ and

\noindent
$(g,1_2)((g',0)(t,s))=(g,1_2)((\theta(1_2)g')t,\mbox{sgn}(g')s)=(\theta(1_2)g\tilde{h}(\theta(1_2)g')(t),-\mbox{sgn}(g)\mbox{sgn}(g')s)$

\noindent
$=((\theta(1_2)g)g'\tilde{h}(t),-\mbox{sgn}(g)\mbox{sgn}(g')s)=((g,1_2)(g',0))(t,s)$;

\vspace{1mm}

\noindent
(iii) $((g,1_2)(g',1_2))(t,s)=(g\theta(1_2)g',0)(t,s)=((\theta(1_2)g\theta(1_2)g')t,-\mbox{sgn}(g\theta(1_2)g')s)=
((\theta(1_2)g)g')t,-\mbox{sgn}(gg')s)$ and

\noindent
$(g,1_2)((g',1_2)(t,s))=(g,1_2)((\theta(1_2)g')\tilde{h}(t),-\mbox{sgn}(g')s)=$

\noindent
$((\theta(1_2)g)\tilde{h}((\theta(1_2)g')\tilde{h}(t)),-\mbox{sgn}(g)\mbox{sgn}(g')s)
=((\theta(1_2)g)g')t,-\mbox{sgn}(g)\mbox{sgn}(g')s)=$

\noindent
$((g,1_2)(g',1_2))(t,s)$.
\par Consequently, $\circ : (F\rtimes_\theta\mathbb{Z}_2)\times(\tilde{\Gamma}\times\mathbb{S}^{2n})
\to\tilde{\Gamma}\times\mathbb{S}^{2n}$ is a well-defined action. Because it does not exist $g\in F$ such that
 $\left\{\begin{array}{l}
\theta(1_2)(g)=g^{-1},\\
\varphi(g,0)=1_2
\end{array}
\right.$
for any $g\in F$, the action  $\circ : (F_m\rtimes_\theta\mathbb{Z}_2)\times(\tilde{\Gamma}\times\mathbb{S}^{2n})
\to\tilde{\Gamma}\times\mathbb{S}^{2n}$ is free. Otherwise suppose that $(g,1_2)\circ(t,s)=(t,s)$. Then we have $(t,s)=((\theta(1_2)g)(\tilde{h}(t)),-\mbox{sgn}(g)s)=\tilde{h}(gt)$ which implies $\mbox{sgn}(g)=-1$ and
$t=(\theta(1_2)g)(\tilde{h}(t))$. The second equation is equivalent to $\tilde h(t)=\tilde{h}^2(gt)=gt=g\theta(1_2)(g)\tilde h(t)$ or $g\theta(1_2)(g)=1$. So the system of equations has a solution which is a contradiction.
So we have  a free, properly discontinuous and cellular action. Further, the induced homomorphism
$\varphi : F\rtimes_\theta\mathbb{Z}_2\to \mbox{Aut}\,(H^{2n}(\tilde{\Gamma}\times\mathbb{S}^{2n}),\mathbb{Z})$
coincides with the given one $\varphi : F\rtimes_\theta\mathbb{Z}_2\to\mbox{Aut}\,(\mathbb{Z})$.

\vspace{3mm}

\par Next, suppose that $\left\{\begin{array}{l}
\theta(1_2)(g)=g^{-1},\\
\varphi(g,0)=1_2
\end{array}
\right.$
for some $g\in F$ and there is an action $(F\rtimes_\theta\mathbb{Z}_2)\times\Sigma(2n)\to\Sigma(2n)$.
Then, on  one hand we have that $\varphi(g,1_2)=\varphi(g,0)\varphi(e,1_2)=0$
and on  the other hand, Proposition \ref{pro}(1) leads to $\varphi(g,1_2)=1_2$,
because the order of $(g,1_2)\in F\rtimes_\theta\mathbb{Z}_2$ is two. This contradiction completes the proof.

\hfill$\square$

\begin{CorX} The group $F\rtimes_\theta\mathbb{Z}_2$ acts on $\Sigma(2n)=\tilde{\Gamma}\times\mathbb{S}^{2n}$ for
any $n\ge 1$, where $\Gamma$ is a graph $($a finite graph provided $F$ is of finite rank$)$ with $\pi_1(\Gamma)=F$.\label{CC}
\end{CorX}
{\bf Proof.} Given the group $F\rtimes_\theta\mathbb{Z}_2$,  consider the homomorphism
$\varphi : F\rtimes_\theta\mathbb{Z}_2\to\mbox{Aut}\,(\mathbb{Z})$ given by the projection map
onto the second factor. Then, in view of Lemma \ref{lem}, the pair $(\theta,\varphi)$ is realizable
and this leads to an action of $F\rtimes_\theta\mathbb{Z}_2$ on
$\Sigma(2n)=\tilde{\Gamma}\times\mathbb{S}^{2n}$ for any $n\ge 1$, where $\Gamma$ is a graph (a finite one provided $F$ is of finite rank)
with $\pi_1(\Gamma)=F$, and the result follows.

\hfill$\square$

\par Now, let $F_m$ be the free group with finite rank $m\ge 1$. For $\theta :\mathbb{Z}_2\to\mbox{Aut}\,(F_m)$ and
$\varphi : F_m\rtimes_\theta\mathbb{Z}_2\to\mbox{Aut}\,(\mathbb{Z})$
with $\varphi|_{\mathbb{Z}_2}=\mbox{id}_{\mathbb{Z}_2}$, we classify
all realizable pairs $(\theta,\varphi)$, i.e., in view of Lemma \ref{lem},  pairs $(\theta,\varphi)$
for which it does not exist  $g\in F_m$ such that
$$\left\{\begin{array}{l}
\theta(1_2)(g)=g^{-1},\\
\varphi(g,0)=1_2.
\end{array}
\right.$$

\par First, we recall a very useful result  by Dyer and Scott \cite[Theorem 3]{DS}.
\begin{thmX}
Let $F$ be any free group, $\theta :\mathbb{Z}_2\to\mbox{\em Aut}(F)$ a homomorphism and
$F^{\theta(1_2)}<F$ the fixed point subgroup of the automorphism $\theta(1_2)$.
Then there is a decomposition \label{ds}
$$F=F^{\theta(1_2)}\ast(\ast_{i\in I}F_i)\ast(\ast_{\lambda\in\Lambda}F_\lambda)$$
into the free product, where each factor is $\theta(1_2)$-invariant and:

\mbox{\em (i)} for each $i\in I$, $F_i=\big<x_{i,1},x_{i,2}\big>$ such that $$\theta(1_2)(x_{i,r})=x_{i,r+1\, (\bmod\,2)}\;\;\mbox{for}\;\;r=1,2;$$

\mbox{\em (ii)} for each $\lambda\in\Lambda$, there is a set $J_\lambda$ with $F_\lambda=\big<x_{\lambda},y_j\mid\,j\in J_\lambda\big>$ such that
$$\theta(1_2)(x_{\lambda})=x_{\lambda}^{-1}\;\mbox{and}\;$$
$$\theta(1_2)(y_j)=x_{\lambda}^{-1}y_jx_{\lambda}\;\mbox{for $j\in J_\lambda$ and $\lambda\in \Lambda$}.$$
\end{thmX}

\bigskip

Basing on Lemma \ref{lem} and Theorem \ref{ds}, we can provide a criterion to decide
whether  a pair is realizable or not. For this purpose  the following is useful.
\par A well-known representation of $\mbox{Aut}\,(F_m)$ is given by
$$\rho_m : \mbox{Aut}\,(F_m)\to \mbox{Aut}\,(F_m/F'_m)\cong GL_m(\mathbb{Z}),$$
where $F'_m$ is the commutator subgroup of $F_m$,
$GL_m(\mathbb{Z})$ the group of all invertible $m\times m$-matrices over $\mathbb{Z}$
and $\rho_m(\theta)$ is the automorphism of the free abelian group $F_m/F'_m\cong\mathbb{Z}^m$ induced by $\theta\in \mbox{Aut}\,(F_m)$.
In view of \cite{Niel}, the group $\mbox{Aut}\,(F_m)$ is finitely presented and $\rho_m$ is surjective.
Because the inner automorphisms $\mbox{Inn}(F_m)\subseteq \mbox{ker}\,\rho_m$,
there is the induced homomorphism $\tilde{\rho}_m : \mbox{Out}(F_m)=\mbox{Aut}\,(F_m)/\mbox{Inn}(F_m)\to GL_m(\mathbb{Z})$
the kernel of which is called the classical {\em Torelli group} denoted by $\mathcal{T}_m$.
\par  Write $I_m$ for the identity $m\times m$-matrix and define $m\times m$-matrices:

$$A(k,r,s)=\begin{pmatrix}
\begin{pmatrix}
0&1\\
1&0\\
\end{pmatrix}&\mathbf{0}&&&&&\\
&&\ddots&&&\\
&&\mathbf{0}&\begin{pmatrix}
0&1\\
1&0\\
\end{pmatrix}&\mathbf{0}&\\
&&&\mathbf{0}&I_r&\mathbf{0}&\\
&&&&\mathbf{0}&-I_s&\\
\end{pmatrix}$$

\vspace{2mm}

\noindent
over integers which satisfy $A(k,r,s)^2=I_m$ with $k$ matrices $\begin{pmatrix}0&1\\1&0\end{pmatrix}$ and $m=2k+r+s$.

Given $\theta : \mathbb{Z}_2\to \mbox{Aut}(F_m)$ with $F_m=\big<x_1,\ldots,x_m\big>$, by  $\rho_m(\theta(1_2))$ we intend to denote the
matrix of the automorphism of the abelianization $F_m^{ab}\cong\mathbb{Z}^m$ with respect to the basis $\{\bar{x}_1,\ldots,\bar{x}_m\}$, where
$\bar{x_i}$ is the projection of $x_i$ onto $\mathbb{Z}^m$ for $i =1,\ldots,m$.
Because $\theta^2(1_2)= \mbox{id}_{F_m}$ we have that $\rho_m(\theta(1_2))^2=I_m$.
\par If $\theta: F_m\to F_m$ is given as follows:
$\theta(x_l)=x_{l+1}$ and $\theta(x_{l+1})=x_l$ for $l=1,3,\ldots,2k-1$;
$\theta(x_l)=x_l$ for $l=2k+1,2k+2,\ldots,2k+r$ and $\theta(x_l)=x_l^{-1}$ for
$l=2k+r+1,2k+r+2,\ldots,2k+r+s$ then  $\rho_m(\theta(1_2))=A(k,r,s)$.
Given $g\in F_m=\big<x_1,\ldots,x_m\big>$, write $|g|_{x_i}$ for its {\em $x_i$-exponent}, i.e.,
$|g|_{x_i}=\sum_{k=1}^tr_k$ provided $x_i^{r_k}$  appears in $g$ for $i=1,\ldots,m$
and $k=1,\ldots,t$, and zero otherwise. Now,
%in virtue of Proposition \ref{ds} and Lemma \ref{r} and making use of the above,
we are in position to state:

\begin{thmX}  Let $F_m=\big<x_1,\ldots,x_m\big>$ be  a free group with $m\ge 1$, $\theta
: \mathbb{Z}_2\to \mbox{\em Aut}(F_m)$ and $\varphi : F_m\rtimes_\theta
\mathbb{Z}_2\to\mbox{\em Aut}(\mathbb{Z})$ be homomorphisms  such that $\rho_m(\theta(1_2))=A(k,r,s)$ and
${\varphi}|_{\mathbb{Z}_2}=\mbox{\em id}_{\mathbb{Z}_2}$. Then the pair $(\theta,\varphi)$ is realizable \label{p}
if and only if $\varphi(x_l,0)=0$ for $l=2k+r+1,\ldots,2k+r+s$.
\end{thmX}
{\bf Proof.} Consider the elements $g\in F_m$ such $\theta(g)=g^{-1}$. If $g$ belongs to the commutator subgroup of $F_m$ then
$\varphi(g,0)=0$. Because $\rho_m(\theta(1_2))=A(k,r,s)$, the equation $\theta(g)=g^{-1}$ implies $A(k,r,s)(\bar g)=-\bar g$,
where $\bar g$ is the projection of $g$ into $\mathbb{Z}^m$. So we get:

\vspace{2mm}

$\begin{pmatrix}
\begin{pmatrix}
0&1\\
1&0\\
\end{pmatrix}&\mathbf{0}&&&&&\\
&&\ddots&&&\\
&&\mathbf{0}&\begin{pmatrix}
0&1\\
1&0\\
\end{pmatrix}&\mathbf{0}&\\
&&&\mathbf{0}&I_r&\mathbf{0}&\\
&&&&\mathbf{0}&-I_s&\\
\end{pmatrix}
\begin{pmatrix}|g|_{x_1}\\\vdots\\|g|_{x_{2k}}\\\vdots\\|g|_{x_{2k+r}}\\\vdots\\|g|_{x_{2k+r+s}}
\end{pmatrix}
=\begin{pmatrix}-|g|_{x_1}\\\vdots\\-|g|_{x_{2k}}\\\vdots\\-|g|_{x_{2k+r}}\\\vdots\\-|g|_{x_{2k+r+s}}
\end{pmatrix}.$

\vspace{2mm}

Consequently, $\left\{\begin{array}{l}
|g|_{x_1}=-|g|_{x_2},\\
\vdots\\
|g|_{x_{2k-1}}=-g|_{x_{2k}}
\end{array}
\right.$ and $\left\{\begin{array}{l}
|g|_{x_{2k+1}}=0,\\
\vdots\\
|g|_{x_{2k+r}}=0.
\end{array}
\right.$

\vspace{1mm}

\par Because $\varphi(x_i,0)=\varphi(\theta(1_2)(x_i),0)=\sum_{l=1}^{2n+r}|\theta(1_2)(x_i)|_{x_l}\varphi(x_l,0)$ for $i=1,\ldots,2k$,
we derive that $\varphi(x_i,0)=
\left\{\begin{array}{ll}
\varphi(x_{i+1},0)&\mbox{if $i$ is odd}\\
\varphi(x_{i-1},0)&\mbox{if $i$ is even}.\\
\end{array}
\right.$

 Now we show one implication. Suppose that  $\varphi(x_l,0)=0$ for $l=2k+r+1,\ldots,2k+r+s$
and let us assume that $g\in F_m$ is a solution of the equation $\theta(g)=g^{-1}$.
Then $\varphi(g,0)=\sum_{i=1}^{2k+r+s}|g|_{x_i}\varphi(x_i,0)=
\sum_{i=1}^{2k}|g|_{x_i}\varphi(x_i,0)=$

\noindent
$|g|_{x_1}\varphi(x_1,0)-|g|_{x_1}\varphi(x_1,0)+
\cdots+|g|_{x_{2k-1}}\varphi(x_{2k-1},0)-|g|_{x_{2k-1}}\varphi(x_{2k-1},0)
=0$.

\noindent
Therefore the system given by the Lemma \ref{lem} has no solution and the result follows.

 To show the converse  suppose  the pair $(\theta,\varphi)$ is realizable. We know that
$x_l$   for $l=2k+r+1,\ldots,2k+r+s$ satisfies the   equation $\theta(g)=g^{-1}$. Therefore by Lemma 3.1 it follows that   $\varphi(x_l,0)=0$ for $l=2k+r+1,\ldots,2k+r+s$ and the proof is complete.

\hfill$\square$

\begin{RemX}{\em (1) It  is not difficult to show the  Theorem \ref{p} for
$F\rtimes_\theta \mathbb{Z}_2$, where $F$ is a free group of arbitrary rank, once we
adpte the hypothesis. Namely we assume that
 $\rho(\theta(1_2))=A(k,r,s)$,  where now we allow that $k,r,s$ can be infinite cardinals,
and  $\varphi(x_l,0)=0$ for all indices $l$ which correspond
to those, where the diagonal is $-1$.}\\
{\em (2) Theorem \ref{ds} tells that there is at least one basis for $F_m$ such that the hypothesis
 $\rho_m(\theta(1_2))=A(k,r,s)$ holds provided $m$ is finite. This is not clear if $F$ is of infinite rank. Although
any free group $F$ and, by Corollary \ref{CC}, the group $F\rtimes_{\theta}\Z_2$ act on a homotopy sphere $\Sigma(2n)$.}
\end{RemX}
\renewcommand{\thesection}{{}}
\section{}
\renewcommand{\thesection}{\arabic{section}}

{\large\bf 4.\ Miscellanea.} A very good survey about the subject below one can  find in \cite{HP}, where several questions are posed and discussed.
In this section, we use some information from \cite{HP} and make related comments having in mind mainly actions on  homotopy spheres
$\Sigma(2n)$.

We begin by recalling that by \cite[Corollary 5.6]{MT}, the Thompson group $$F=\big<x_0,x_1,\ldots|\,x_ix_jx_i^{-1}=x_{j+1},\,i<j\big>$$
with $\mbox{cd}\,F=\mbox{vcd}\,F=\infty$ does not act freely and properly discontinuously on any $\mathbb{R}^m\times\mathbb{S}^n$.
In fact, it is also  true that $F$ does not act freely and properly discontinuously on any homotopy sphere $\Sigma(n)$.
To see this, suppose that $F$ acts on some $\Sigma(n)$. Using the fact that for every positive integer $m$ the group $F$
contains a copy of the free abelian group $\Z^m$, Proposition \ref{PP} implies the inequality $\dim \Sigma(n)\geq \mbox{vcd}\,\Z^m+n=n+m$.
But this is not possible for an arbitrary $m$ and the result follows.
\par In view of \cite{KK}, this countable group $F$ has periodic cohomology in the sense that $H^i(F, \Z) \cong H^{i+2}(F,\Z)$ for all $i>1$.
But we are unable to show its periodicity after some steps (in the sense of \cite{T1}). Using \cite{KK}, it can be shown,
without too much difficulties that this isomorphism cannot be realized by means of the cup product.

\par S.\ Prassidis has shown in the paragraph following  \cite[Theorem 10]{PS} that:
\begin{thmX}There exist discrete groups $G$ with $\mbox{\em vcd}\,G=\infty$ which act freely and properly on some $\mathbb{R}^m\times\mathbb{S}^n$.
\end{thmX}

The action given in \cite{PS} is free, properly discontinuous but not  co-compact. Then  F.T.\ Farrell and C.W.\ Stark \cite[Theorem 1]{FS}
showed:
\begin{thmX} For each $m\ge 2$ and $n\ge m(m+1)$, there are smooth closed manifolds with universal covering spaces
$\mathbb{R}^m\times\mathbb{S}^{2n-1}$ and fundamental group of infinite virtual cohomological dimension.
\end{thmX}

Groups from the results above are torsion with $\mbox{vcd}\,G=\infty$ and in view of Proposition \ref{pr} they cannot act on any $\Sigma(2n)$, in particular
on any $\R^m\times\S^{2n}$.
But it is natural to ask: can a  torsion-free group  $G$ with $\mbox{cd}\,G=\infty$ acts (possibly co-compactly),  freely and  properly discontinuously
on some $\mathbb{R}^m\times\mathbb{S}^n$? By private communication with  F.X.\ Connolly and S.\ Prasidis  this  question is unsettled.

Several  of the questions and results above, can be studied  if we restrict to the
family of homotopy spheres $\Sigma(2n)$. Taking into account \cite[Theorem 5.2]{L}, we close this paper with:
\begin{questX} {\em Suppose that a group $G$ acts, freely and properly discontinuously
(possibly co-compactly) on some $\Sigma(2n)$ with $\mbox{dim}\,\Sigma(2n)\le m+2n$.\label{q}
Does it follow that $\mbox{vcd}\,G\le m$?}
\end{questX}
Certainly, the proof of \cite[Theorem 5.2 (1)]{L} leads to $\mbox{vcd}\,G\le m$ provided $\mbox{vcd}\,G<\infty$. Notice that Proposition \ref{Pp} yields
$\mbox{vcd}\,G\le 1$ for any virtually cyclic group $G$ and $\mbox{dim}\,\Sigma(2n)<\infty$. Further, by the proof of \cite[Corollary 7.2]{L}
the answer to Question \ref{q} is affirmative for $m=0,1$.

{}

\vspace{2cm}

\noindent{\bf Institute of Mathematics\\
Casimir the Great University\\
pl.\ Weyssenhoffa 11\\
85-072 Bydgoszcz, Poland}\\
e-mail: marek@ukw.edu.pl

\vspace{1.5mm}

\noindent{\bf Department of Mathematics-IME\\
University of S\~ao Paulo\\
Caixa Postal 66.281-AG. Cidade de S\~ao Paulo\\
05314-970 S\~ao Paulo, Brasil}\\
e-mail: dlgoncal@ime.usp.br

\vspace{1.5mm}

\noindent
{\bf Instituto de Matem\'aticas, Unidad Oaxaca\\
Universidad Nacional Aut\'onoma de M\'exico\\
 Oaxaca, Oax.\ M\'exico}\\
e-mail: rolando@matcuer.unam.mx

\end{document}